\documentclass{article}
\usepackage{amsmath, amssymb}
\usepackage{amssymb,pb-diagram,pst-all}
\usepackage{amscd}

\usepackage{fancybox}
\newcommand{\ZZ}{\mathbb Z}
\newcommand{\PP}{\mathbb P}
\newcommand{\QQ}{\mathbb Q}
\newcommand{\CC}{\mathbb C}

\newcommand{\mcB}{\mathcal B}

\newcommand{\mcD}{\mathcal D}
\newcommand{\mcQ}{\mathcal Q}

\newcommand{\tor}{{\mathrm {tor}}}

\newcommand{\lex}{{\mathrm {lex}}}

\newcommand{\Cov}{\mathop {\rm Cov}\nolimits}

\newcommand{\MW}{\mathop {\rm MW}\nolimits}
\newcommand{\Gal}{\mathop {\rm Gal}\nolimits}

\newcommand{\NS}{\mathop {\rm NS}\nolimits}

\newcommand{\rank}{\mathop {\rm rank}\nolimits}
\newcommand{\Red}{\mathop {\rm Red}\nolimits}
\newcommand{\Pic}{\mathop {\rm Pic}\nolimits}

\newcommand{\Supp}{\mathop {\rm Supp}\nolimits}
\newcommand{\Sing}{\mathop {\rm Sing}\nolimits}

\def\miya#1{\mathrel{\mathop{\rightarrow}\limits^#1}}

\newcommand{\mcC}{\mathcal C}

\newtheorem{thm}{Theorem}[section]
\newtheorem{cor}[thm]{Corollary}
\newtheorem{prop}[thm]{Proposition}
\newtheorem{lem}[thm]{Lemma}

\newtheorem{exmple}[thm]{Example}
\newtheorem{rem}[thm]{Remark}
\newtheorem{qz}[thm]{Question}

\newtheorem{mthm}{Theorem}

\newtheorem{Rem}[mthm]{Remark}

\newtheorem{Definition}[mthm]{Definition}

\newcommand{\I}{\mathop {\rm I}\nolimits}

\newcommand{\qed}{\hfill $\Box$}

\newcommand{\proof}{\noindent{\textsl {Proof}.}\hskip 3pt}

\renewcommand{\thesubparagraph}{\theparagraph.\@arabic\c@subparagraph}
\begin{document}
  
  \begin{center}
  
 {\bf  \Large 
 On the Abel-Jacobi map for bisections of \\ rational elliptic surfaces \\
 and \\
 Zariski $N$-plets for conic arrangements}
\bigskip

\bigskip
\large 
Shinzo BANNAI 
and 
Hiro-o TOKUNAGA\footnote{Research partly supported by the research grant 22540052
from JSPS.  MSC 2010 Classification: 14E20, 14J27}

\end{center}
\normalsize

\begin{abstract}
We study the Abel-Jacobi map for bisections of a certain rational elliptic surface. As an application, we construct examples of Zariski $N$-plets for conic arrangements.
\end{abstract}

\bigskip

\section*{Introduction}

In this article, we study the geometry of bisections and the Able-Jacobi map for them on rational elliptic
 surfaces. We apply our result to construct certain conic arrangements and to
 study the topology of the conic arrangement via Galois covers. As an application,
 we give  examples of 
Zariski $N$-plets. 

Let $\varphi : S \to C$ be an elliptic surface over a smooth projective curve $C$ with
a section $O$ such that no $(-1)$ curve is contained in any fiber and $\varphi$ has
at least one degenerate fiber (see \S 1 for terminologies on elliptic surfaces). We denote
the set of sections of $\varphi$ by $\MW(S)$. The following facts are well-known:

\begin{enumerate}

\item[(I)] By regarding $O$ as the zero element and considering fiberwise addition, $\MW(S)$ becomes
an abelian group. We denote its addition and the multiplication-by-$m$-map by
$\dot{+}$ and $[m]$, respectively.

\item[(II)] Abel's theorem. Let $D$ be a divisor on $S$ and put $d =DF$, where $F$ is
a fiber of $\varphi$. There exists a unique section $s(D)$ (the image of the Abel-Jacobi map)
 such that $D - s(D) - (d-1)O$ is
algebraically equivalent to a divisor whose irreducible components are all in fibers of $\varphi$.
\end{enumerate}

From these facts, by identifying a section and its image, we immediately observe the following:

\begin{enumerate}

\item[(a)] Given $s_1, s_2 \in \MW(S)$, then we have new curves $s_1\dot{+}s_2$ and 
$[m]s_i$ 
$(i = 1, 2)$ on $S$.

\item[(b)] Given a divisor $D$, then we have a new curve $s(D)$ on $S$.

\end{enumerate} 

In \cite{tokunaga12}, we make use of observation (a) in order to study Zariski pairs for
conic-line arrangements.  In this article, however, we study the 
relation between $D$ and $s(D)$
in observation (b) in the case where $S$ is a rational elliptic surface with an 
$\I_2$ fiber  and $D$ is a bi-section
in order to give conic arrangements with prescribed conditions. 

We now explain our main result.
Let $\mcQ$ be a reduced plane quartic curve. $\mcQ$ can be reducible, but we will assume that
$\mcQ$ has a component of degree greater than or equal to $2$.  Let $S^{\prime\prime}$  be
the double cover of $\PP^2$ branched 
along $\mcQ$ and $\bar{S}$ the canonical resolution of the singularities of $S^{\prime\prime}$ 
(see \cite{horikawa} for the canonical resolution). Let $p_1$ and $p_2$ be  
general points on a component of $\mcQ$ with degree $\ge 2$. We will denote the inverse images of $p_1$ and $p_2$ in $\bar{S}$ by the same symbol. 
The pencils  of lines through $p_i$  $(i=1,2)$ give rise to  pencils of elliptic curves through 
$p_i$ on $\bar{S}$. By 
resolving the base points of each pencil we obtain  rational elliptic surfaces
$\varphi_{p_1} : S_{p_1} \to \PP^1$ and $\varphi_{p_2} : S_{p_2} \to \PP^1$. 
The resolution maps will be denoted by $\mu_i:S_{p_i}\rightarrow \bar{S}$. Each $\mu_i$ is a composition of
 two blowing ups. 

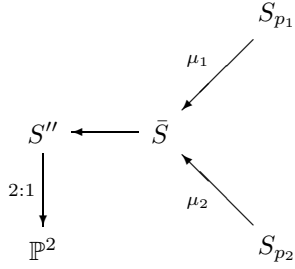
\[
\begin{diagram}
\node{}\node{}\node{S_{p_1}}\arrow{sw,t}{\mu_1}\\
\node{S^{\prime\prime}}\arrow[1]{s,l}{2:1}\node{\bar{S}}\arrow{w}\node{}\\
\node{\PP^2}\node{}\node{S_{p_2}}\arrow{nw,b}{\mu_2}
\end{diagram}
\]

The exceptional divisor of the second blowing-up of $\mu_1$ (resp. $\mu_2$) gives rise to a section of 
$S_{p_1}$ (resp. $S_{p_2}$). We will regard this section as the zero section and denote it by 
$O_{1}$ (resp. $O_2$). By construction, $S_{p_1}$ and $S_{p_2}$ are rational elliptic surfaces that have the 
same configuration of singular fibers. 
%Each surface has a $I_2$ type singular fiber whose components arise from the exceptional divisor of the first blowing-up in $\mu_i$ which meets $O_{i}$ and the strict transform
% of the tangent line $T_{p_i}$ of $Q$ at $p_i$. We will denote these components by $\Theta_{p_i,0}$ and 
% $\Theta_{p_i,1}$.  
% 
% All the other reducible singular fibers arise from the exceptional sets of the resolution 
% $S\rightarrow S^\prime$, hence they are in 1 to 1 correspondence with the singularities of $Q$. We will
%  denote their components by $\Theta_{v,i}$ where $v\in \Sing(Q)$.

Let $D_1,\ldots,D_m$ be divisors on $\bar{S}$ such that they do not pass through  $p_1, p_2$ and 
their strict transforms under $\mu_1$ (resp.  $\mu_2$) give rise to sections of 
$S_{p_1}$ (resp. $S_{p_2}$). Let $s_i(D_j)$ denote the section corresponding to $D_j$ on $S_{p_i}$. 
Let $C_i=s_i(D_1)\dot+\cdots\dot+s_i(D_m)$. 
Put $\overline{C}_i = \mu_i(C_i)$ and let $\widehat{C}_i$ be the strict 
transform of $\overline{C}_i$ under $\mu_j^{-1}$.  $(i\not=j)$. For general $p_1$ and $p_2$,
 $\widehat{C}_2$ becomes a multi-section of $S_{p_1}$. 
 Under this setting, we have:

\begin{mthm}\label{prescription}
{
Suppose that $C_2\not=O_2$ and $p_1\not\in \overline{C}_2$. Then
\[
s(\widehat{C}_2)=C_1.
\]
%where $\bar{\psi}_1$ is the map from $\NS(S_{p_1})$ to $\MW(S)$ given in \cite{shioda90}.
}
\end{mthm}

We apply Theorem~\ref{prescription} to study dihedral covers of $\PP^2$ whose 
branch locus is a conic arrangement and give some examples  of Zariski $N$-plets
for conic arrangements.

Let us first recall the definition of a Zariski $N$-plet.

\begin{Definition}\label{def:zariski-N-plet}{ \rm An $N$-plet of reduced plane curves
$(\mcB_1, \ldots, \mcB_N)$ in $\PP^2$ is said to be a Zariski $N$-plet if it satisfies the
following conditions:

\begin{enumerate}

\item[(i)] For each $i$, there exists a tubular neighborhood $T(\mcB_i)$ of $\mcB_i$ such
that $(T(\mcB_i), \mcB_i)$ is homeomorphic to $(T(\mcB_j), \mcB_j)$ for any 
$1 \le i < j \le N$.

\item[(ii)] For any $1 \le i < j \le N$, there exists no homeomorphism from 
$(\PP^2, \mcB_i)$ to $(\PP^2, \mcB_j)$.
\end{enumerate}

}
\end{Definition}

The first condition can be replaced by the combinatorial type of each curve. If 
we denote the irreducible decomposition of $\mcB_i$ by $\mcB_i = \mcC_{i,1} + 
\ldots + \mcC_{i, r_i}$, the combinatorial type of $\mcB_i$ is, roughly speaking, determined
by $\deg \mcC_{i,j}$, the set of topological types of the 
singularities of $\mcC_{i, j}$ and how
the irreducible components meet each other (For details, see \cite{act}).
 
 The combinatorial type of the conic arrangement which we consider in this article is as follows:
 
 \begin{Definition}\label{def:k-NT}{\rm A reduced plane curve $\mcB$ consisting of $(k+2)$ 
 irreducible conics is called a Namba-Tsuchihashi conic arrangement of type $k$ ($k$-NT
 arrangement, for short) if  $\mcB$ is of the form $\mcQ + \mcC_1 + \ldots + \mcC_k$ 
 satisfying the following conditions:
  
 \begin{enumerate}
 
  \item[(i)] $\mcQ$ is a quartic consisting of $2$ irreducible conics $\mcC^\prime$, $\mcC^{\prime\prime}$ intersecting transversely.

  \item[(ii)] $\mcC^\prime$, $\mcC^{\prime\prime}$ are tangent to $\mcC_j$ $(j = 1, \ldots, k)$.
 
  \item[(iii)]  $\mcC_i$ and $\mcC_j$ $(1 \le i < j \le k)$ intersect transversally.
  
  \item[(iv)] The singularities of $\mcB$ are only nodes and tacnodes, i.e, 
  each component of $\mcQ$ is tangent
  to $\mcC_j$ ($j = 1, \ldots, k$) at two distinct  points and no three conics meet at one point.
  
    \end{enumerate}
   For $\mcB$ satisfying only the 
   first two conditions, we call $\mcB$ {\it a weak Namba-Tsuchihashi
   arrangement of type $k$} (a weak $k$-NT arrangement, for short).
 }
 \end{Definition}

In \cite{namba-tsuchi}, Namba and Tsuchihashi give an example of a Zariski pair
for Namba-Tsuchihasi arrangements of type $2$. In this article, as an application of
 Theorem~\ref{prescription} and  \cite[Theorems~3.1 and ~4.1]{tokunaga12})
we generalize Namba-Tsuchihashi's example and prove the following:

\begin{mthm}\label{thm:zariski-kplet}{ Let $y(k,3)$ be the number of Young diagrams with
$k$ boxes and at most $3$ rows. There exist $y(k, 3)$ $k$-NT arrangements
$\mcB_1, \ldots, \mcB_{y(k,3)}$ such that no homeomorphism 
$h : (\PP^2, \mcB_i) \to (\PP^2, \mcB_j)$ with $h(\mcQ) = \mcQ$ exists for
any $i < j$. In particular, they form a Zariski $y(k, 3)$-plet for $k \ge 3$.
} 
\end{mthm}

\begin{Rem} {\rm Professor Miles Reid told the second author that $y(k, 3)$ can be computed by the orbifold Riemann-Roch
formula. In fact,
since $y(k, 3)$ satisfies the recursive formula
\[
y(k, 3) = 1 + \left [\frac k2\right ] + y(k-3, 3)
\] for $k \ge 4$, where $[\bullet]$ denotes the maximum integer not exceeding $\bullet$,
one can compute explicit formulas for $y(k, 3)$ as follows:
\[
y(k, 3) = \left \{ \begin{array}{ll}
                         \frac 1{12}(k^2 + 6k + 12) & k \equiv 0 \bmod 6 \\
                         \frac 1{12}(k +1)(k +  5) & k \equiv \pm 1 \bmod 6 \\
                         \frac 1{12}(k + 2)(k + 4 ) & k \equiv \pm 2 \bmod 6 \\
                         \frac 1{12}(k+ 3)^2      & k \equiv 3 \bmod 6
                         \end{array}
                         \right .
\]
The authors thank Professor Miles Reid for his comments.
}
\end{Rem}
                          
This article consists of $4$ sections. We summarize some necessary facts on 
the theory of 
elliptic surfaces and dihedral covers in $\S 1$ and prove Theorem~\ref{prescription}
in \S 2.
In \S 3, we construct $y(k, 3)$ weak $k$-NT arrangement and prove Theorem~\ref{thm:zariski-kplet}
in \S 4.

 %%%%%%%%%%%%%%%%%%%%%%%%%%%%%%%%%%
 
 %Preliminaries

\section{Preliminaries}

\subsection{Elliptic sufaces}
\subsubsection{General Facts}
We first summarize some facts from the theory of elliptic surfaces. As for details, we refer to 
\cite{kodaira}, \cite{miranda-basic}, \cite{miranda-persson} and \cite{shioda90}.

 In this article, by an {\it elliptic surface}, we always mean a smooth projective surface $S$ 
  with a fibration $\varphi : S \to C$ over a smooth projective curve, $C$, as follows:
 \begin{enumerate}
  \item[(i)] There exists a non-empty finite subset, $\Sing(\varphi)$, of $C$ such
  that $\varphi^{-1}(v)$ is a smooth curve of genus $1$ (resp. a singular curve )for $v \in C\setminus \Sing(\varphi)$ (resp. $v \in \Sing(\varphi))$
 \item[(ii)]  $\varphi$ has a section
 $O : C \to S$ (we identify $O$ with its image). 
 \item[(iii)]  $\varphi$ is minimal, i.e., there is no exceptional
 curve of the first kind in any fiber. 
 \end{enumerate}
 
 For $v \in \Sing(\varphi)$, we put $F_v = \varphi^{-1}(v)$. 
 We denote its irreducible decomposition by 
 \[
 F_v = \Theta_{v, 0} + \sum_{i=1}^{m_v-1} a_{v,i}\Theta_{v,i}, 
 \]
 where $m_v$ is the number of irreducible components of $F_v$ and $\Theta_{v,0}$ is the
 irreducible component with $\Theta_{v,0}O = 1$. We call $\Theta_{v,0}$ the {\it identity
 component}.  The classification  of singular fibers is well known(\cite{kodaira}). %The types we will use in this paper are as follows:
 Note that every
 smooth irreducible component of reducible singular fibers is a rational curve with self-intersection number $-2$.
 
% \begin{center}
 
 %%%%%%%%%%%%%%%%%%%%%%%%
%\input{Ib-fiber.tex}
%%%%%%%%%%%%%%%%%%%%%%%%%
 
% %%%%%%%%%%%%%%%%%%%%%%%%
%\input singularfiber2.tex
%%%%%%%%%%%%%%%%%%%%%%%%%%

% 
% %%%%%%%%%%%%%%%%%%%%%%%%
%\input singularfiber3.tex
%%%%%%%%%%%%%%%%%%%%%%%%%%

%
%\end{center} 

 We also define a subset of $\Sing(\varphi)$ by
 $\Red(\varphi) := \{v \in \Sing(\varphi)\mid \mbox{$F_v$ is reducible}\}$. Let 
 $\MW(S)$ be the set of sections of $\varphi : S \to  C$. From our assumption, 
 $\MW(S) \neq  \emptyset$. By regarding $O$ as the zero element of 
 $\MW(S)$ and 
 considering fiberwise addition (see \cite[\S 9]{kodaira} or \cite[\S 1]{tokunaga11} for
 the addition on singular fibers), 
 $\MW(S)$ becomes an abelian group. We denote its addition by $\dot{+}$.
 Note that the ordinary  $+$ is used for the sum of divisors, and the two operations should not be confused. 
 
% Also for $k \in \ZZ$ and $s \in \MW(S)$, we write
% \[
% [m]s := \left \{ \begin{array}{l}
%                  \mbox{$m$-times addition of $s$ if $m \geq 0$} \\
%                  \mbox{$m$-times addition of the inverse of $s$ if $m < 0$}.
%                  \end{array} \right .
%\]
Let $\NS(S)$ be the N\'eron-Severi group of $S$ and let $T_{\varphi}$ be the 
subgroup of $\NS(S)$ generated by $O, F$ and  $\Theta_{v,i}$ $(v \in \Red(\varphi)$,
$1 \le i \le m_v-1)$. Then we have the following theorems:

\begin{thm}\label{thm:shioda-basic0}{\cite[Theorem~1.2]{shioda90}) Under our assumption,
$\NS(S)$ is torsion free.
}
\end{thm}

\begin{thm}\label{thm:shioda-basic}{(\cite[Theorem~1.3]{shioda90}) Under our assumption,
there is a natural map $\tilde{\psi} : \NS(S) \to \MW(S)$ which induces an isomorphisms of 
groups
\[
\psi : \NS(S)/T_{\varphi} \cong \MW(S).
\]
In particular, $\MW(S)$ is a finitely generated abelian group.
}
\end{thm}

In the following,  by the rank of $\MW(S)$, denoted by $\rank\MW(S)$, we mean that of the free part 
of $\MW(S)$.
 For a divisor on $S$, we  put $s(D) = \psi(D)$. Then we have

\begin{lem}\label{lem:fund-relation}{(\cite[Lemma~5.1]{shioda90})
$D$ is uniquely written in the form:
\[
D \approx s(D) + (d-1)O + nF + \sum_{v\in \Red(\varphi)}\sum_{i=1}^{m_v-1}b_{v,i}\Theta_{v,i},
\]
where $\approx$ denotes the algebraic equivalence of divisors, and $d, n$ and $b_{v,i}$
are integers defined as follows:
\[
d = DF \qquad n = (d-1)\chi({\mathcal O}_S) + OD - s(D)O,
\]
and
\[
\left ( \begin{array}{c}
          b_{v,1} \\
          \vdots \\
          b_{v, m_v-1} \end{array} \right ) = A_v^{-1}\left (\begin{array}{c}
                                                                          D\Theta_{v,1} - s(D)\Theta_{v,1} \\
                                                                          \vdots \\
                                                                      D\Theta_{v,m_v-1} - s(D)\Theta_{v, m_v-1}
                                                                      \end{array} \right )
\]
Here $A_v$ is the intersection matrix $(\Theta_{v,i}\Theta_{v, j})_{1\le i, j \le m_v-1}$.  
}
\end{lem}

For a proof, see \cite{shioda90}.

\medskip

%Put $\NS_{\QQ}:= \NS(S)\otimes \QQ$ and $T_{\varphi, \QQ} := T_{\varphi}\otimes \QQ$.
%Since $\NS(S)$ is torsion free under our setting, there is no harm in considering $\NS_{\QQ}$.
%By using the intersection pairing, we have the orthogonal decomposition $\NS_{\QQ}
%= T_{\varphi, \QQ}\oplus(T_{\varphi, \QQ})^{\perp}$.   In \cite{shioda90}, the homomorphism
%$\phi: \MW(S) \to (T_{\varphi, \QQ})^{\perp}\subset\NS_{\QQ}$ is defined as follows:
%\begin{eqnarray*}
%\phi : \MW(S) \ni s &\mapsto & s - O - (sO + \chi({\mathcal O}_S))F \\
%&&  +\sum_{v \in \Red(\varphi)}
%(\Theta_{v,1}, \ldots, \Theta_{v, m_v-1})(-A_v)^{-1}
%\left ( \begin{array}{c}
%        s\Theta_{v,1} \\
%        \vdots \\
%        s\Theta_{v, m_v-1}
%        \end{array} \right ) \in (T_{\varphi, \QQ})^{\perp}.
%\end{eqnarray*}
 Also,  in \cite{shioda90},  a $\QQ$-valued bilinear form $\langle \, , \, \rangle$ on
 $\MW(S)$  is defined by using the intersection pairing on $\NS$.   Here are two basic properties of $\langle \, , \, \rangle$:
\begin{itemize}

\item $\langle s, \, s \rangle \ge 0$ for $\forall s \in \MW(S)$ and the equality holds if and 
only if $s$ is an element of finite order in $\MW(S)$. 

\item An explicit formula for $\langle s_1,
s_2\rangle$ ($s_1, s_2 \in \MW(S)$) is given as follows:
\[
\langle s_1, s_2 \rangle = \chi({\mathcal O}_S) + s_1O + s_2O - s_1s_2 - \sum_{v \in \Red(\varphi)}
\mbox{Contr}_v(s_1, s_2),
\]
where $\mbox{Contr}_v(s_1, s_2)$ is given by
\[
\mbox{Contr}_v(s_1, s_2) = (s_1\Theta_{v,1}, \ldots, s_1\Theta_{v, m_v-1})(-A_v)^{-1}
\left ( \begin{array}{c}
        s_2\Theta_{v,1} \\
        \vdots \\
        s_2\Theta_{v, m_v-1}
        \end{array} \right ).
\]
 As for explicit values of 
$\mbox{Contr}_v(s_1, s_2)$, we refer to \cite[(8.16)]{shioda90}.
\item Let $\MW(S)^0$ be the subgroup of $\MW(S)$ given by 
\[
\MW(S)^0 := \{ s \in \MW(S) \mid \mbox{$s$ meets $\Theta_{v, 0}$ for $\forall v \in \Red(\varphi)$.}\}
\]
$\MW(S)^0$ is called the {\it narrow part} of $\MW(S)$. By the explicit formula as above,
$(\MW(S)^0, \, \langle \, , \, \rangle)$ is a positive definite even integral lattice.
\end{itemize}
      
\subsubsection{Double cover construction of an elliptic surface}\label{subsubsec:double-cover}

We refer to \cite[Lectures III and IV]{miranda-basic} for details.
Let $\varphi: S \to C$ be an elliptic surface. We can regard the generic fiber $S_{\eta}$
of $\varphi$ as  an elliptic curve over $\CC(C)$, the rational function field of $C$, 
under our assumption.
The inverse morphism with respect to
the group law on $S_{\eta}$  induces an involution $[-1]_{\varphi}$ on $S$. 
Let $S/\langle [-1]_{\varphi}
\rangle$ be the quotient by $[-1]_{\varphi}$.  It is known that  $S/\langle[-1]_{\varphi}\rangle$ 
is smooth and we can  blow down  $S/\langle [-1]_{\varphi} \rangle$ to its relatively minimal
model $W$ over $C$ in the following way:

Let us denote
      \begin{itemize}
       \item $f: S \to S/\langle [-1]_{\varphi}\rangle$: the quotient morphism, 
       \item $q^\prime: S/\langle [-1]_{\varphi}\rangle \to W$: the blow down, and 
      \item $S \miya{{\mu^\prime}} S' \miya{{f'}} W$: the Stein factorization of $q^\prime\circ f$. 
      \end{itemize}
      Then we have:

\begin{enumerate}

 \item  The branch locus $\Delta_{f'}$ of $f'$ consists of a section $\Delta_0$ and a
 trisection $T$ such that its singularities are at most simple singularities (see \cite[Chapter II, \S 8]{bpv} for simple singularities
 and their notation) and $\Delta_0\cap T = \emptyset$. 
 
 \item $\Delta_0 + T$ is $2$-divisible in $\Pic(W)$.
 
 \item     
   The morphism  $\mu$ is obtained by contracting all the irreducible components of 
 singular fibers not meeting $O$.

 \end{enumerate}

Conversely, if $\Delta_0$ and $T$ on $W$ satisfy the above conditions, we
obtain an elliptic surface $\varphi : S \to C$, as the canonical%\PP^1'ðC'ɏC³
resolution of a double cover $f' : S' \to W$ with $\Delta_{f'} = \Delta_0
+ T$, and the diagram (see \cite{horikawa} for the canonical resolution):
\[
\begin{CD}
S' @<{\mu^\prime}<< S \\
@V{f'}VV                 @VV{f}V \\
W@<<{q^\prime}< \widehat{W}.
\end{CD}
\]
Here $q$ is a composition of blowing-ups so that $\widehat{W}
= S/\langle [-1]_{\varphi}\rangle$. 
Hence any elliptic surface  is obtained in this way.
In the following, we call the diagram above
{\it the double cover diagram for $S$}.

In the case when $S$ is a rational elliptic surface, $W$ is the Hirzebruch surface, 
 $\Sigma_2$, of degree 
 $d = 2$ and $\Delta_{f'}$ is of the form $\Delta_0 + T$, where
 $\Delta_0$ is a section with $\Delta_0^2 = -2$ and $T \sim 3(\Delta_0 + 2{\mathfrak f})$, 
 ${\mathfrak f}$ being a fiber of the ruling $\Sigma_2 \to {\mathbb P}^1$.
\begin{rem}\label{rem:double-cover-construction}{\rm

\begin{itemize}
  \item For each $v \in \Sing(\varphi)$, the type of $\varphi^{-1}(v)$ is determined by the type of singularity of $T$
  on ${\mathfrak f}_v$ and the relative position between ${\mathfrak f}_v$ and $T$ (see \cite[Table 6.2]{miranda-persson}).
  
  \item Let 
 $\sigma_f$ be the covering transformation of $f$. We can check that $\sigma_f$
  coincides with  $[-1]_{\varphi}$ and how $[-1]_{\varphi}$ acts on irreducible components
  of singular fibers (see \cite[Remark 3.1 (i)]{tokunaga12}).

 % We here
%   summarize 
%   the action 
%  of $\sigma_f$ on irreducible components of singular fibers, which we need later:
%  
%  \begin{center}
%  \begin{tabular}{|c|c|} \hline
%  Type of a singular fiber & The action on irreducible component \\ \hline
%  $\I_b$ & $\displaystyle{\begin{array}{lc}
%                                     \Theta_0 \mapsto \Theta_0  & \\
%                                     \Theta_i \mapsto \Theta_{b-i} & i = 1, \ldots, b-1 
%                                     \end{array}}$ \\
%                                                     
%  $\III$ & $\displaystyle{\begin{array}{cc}
%                                     \Theta_i \mapsto \Theta_i  &  \forall i \\
%                                     \end{array}}$ \\ \hline             
%  \end{tabular}   
%  \end{center}                                
 \end{itemize}
  }
  \end{rem}
  
\subsection{Dihedral Covers.}

\subsubsection{Branched Galois covers}

We first explain our terminology for Galois covers. For normal projective varieties $X$ and $Y$ with
finite morphism $\pi : X \to Y$, we say that $X$ is a Galois cover of $Y$ if the induced
field extension $\CC(X)/\CC(Y)$ is Galois, where $\CC(\bullet)$ denotes the rational function 
field of $Y$.
Recall that the Galois group $\Gal(\CC(X)/\CC(Y))$ acts on $X$ such that $Y$ is obtained as the quotient
space with respect to this action ({\it cf.} \cite[\S 1]{tokunaga97}).
 If the Galois group 
$\Gal(\CC(X)/\CC(Y))$ is isomorphic to a finite group $G$,  we simply call $X$ a $G$-cover of $Y$.
The branch locus of $\pi : X \to Y$, which is denoted by $\Delta_{\pi}$ or $\Delta(X/Y)$,
is a subset of $Y$ consisting of points $y$ of $Y$, over which $\pi$ is not locally
isomorphic. It is well-known that $\Delta_{\pi}$ is an algebraic subset of pure codimension $1$
if $Y$ is smooth (\cite{zariski}).  
%For a reduced divisor on $B$, we identify it with its support

 Now assume that $Y$ is smooth. Let $B$ be a reduced divisor on $Y$ and we denote its irreducible
 decomposition by $B = \sum_{i=1}^rB_i$. We say that a $G$-cover $\pi : X \to Y$ is branched
 at $\sum_{i=1}^r e_iB_i$ if $(i)$ $\Delta_{\pi} = B$ (here we identify $B$ with its support) and $(ii)$ the ramification index
 along $B_i$ is $e_i$ for each $i$, where  the ramification index means the
 one along the smooth part of $B_i$ for each $i$.  Note that the study of $G$-covers is related
 to that of the topology of the complement of $B$, as the proposition below holds: 
 
 \begin{prop}\label{prop:fund}{(\cite[Proposition 3.6]{act}) 
Under the notation as above, let $\gamma_i$ be a meridian around $B_i$, and
$[\gamma_i]$ denote its class in the topological fundamental group 
$\pi_1(Y\setminus B, p_o)$.
If there exists a $G$-cover $\pi : X \to Y$ branched at $e_1B_1 + \cdots + e_rB_r$, 
then there exists a normal subgroup $H_{\pi}$ of $\pi_1(Y\setminus B, p_o)$ 
such that:

\begin{enumerate}
\item[(i)] $[\gamma_i]^{e_i} \in H_{\pi},  [\gamma_i]^k \not\in H_{\pi}, (1 \le k \le e_i -1)$, and
\item[(ii)]  $\pi_1(Y\setminus B, p_o)/H_{\pi} \cong G$.
\end{enumerate}

Conversely, if there exists a normal subgroup $H$ of $\pi_1(Y\setminus B, p_o)$ satisfying the above two conditions for $H_{\pi}$, then there exists a $G$-cover
$\pi_H : X_H \to Y$ branched at $e_1B_1 + \cdots + e_rB_r$.
}
\end{prop}

For $G$-covers $\pi_1 : X_1 \to Y$ and $\pi_2 : X_2 \to Y$, we identify them if there exists an isomorphism
$\Phi : X_1 \to X_2$ such that $\pi_1 = \Phi\circ\pi_2$.  
Under the same notation as in Proposition~\ref{prop:fund}, Put

\begin{itemize}
\item $\Cov (Y, B, G) :=$ the set of isomorphism classes of $G$-covers $\pi : X \to Y$ such that 
$\Delta_{\pi} \subseteq B$, and

\item $\Cov_b(Y, e_{i_1}B_{i_1}+\ldots + e_{i_s}B_{i_s}, G)  :=$
the set of isomorphism classes of $G$-covers $\pi : X \to Y$ branched at $\sum_je_{i_j}B_{i_j}$. Here we
only assume that $\Supp(\sum_je_{i_j}B_{i_j}) \subseteq B$.  
\end{itemize}

Note that if $\Supp\left(\sum_je_{i_j}B_{i_j}\right)\subset B$ then  $\Cov_b(Y, e_{i_1}B_{i_1}+\ldots + e_{i_s}B_{i_s}, G)$ is a subset of $\Cov(Y, B, G)$. 
%Note that $\Cov_b(Y, e_{i_1}B_{i_1}+\ldots + e_{i_s}B_{i_s}, G)$ is a subset of $\Cov(Y, B, G)$. 
By Proposition~\ref{prop:fund}, we have the following:

\begin{prop}\label{prop:1-1}{We keep the notation  in Proposition~\ref{prop:fund}. 
Let $B_i$ ($i =1, 2$) be reduced divisors on $Y$. Let their irreducible decompositions
be denoted by  $B_{i,1} + \ldots + B_{i, r_i}$ ($i = 1, 2$). If there exists a homeomorphism $h$ from $(Y, B_1)$
to $(Y, B_2)$ such that $h(B_{1, j}) = B_{2, j}$ $(j = 1, \ldots, r_1 (= r_2))$, then there exists
a one-to-one correspondence between $\Cov_b(Y, e_{i_1}B_{1, i_1}+\ldots + e_{i_s}B_{1,i_s}, G)$ and
$\Cov_b(Y, e_{i_1}B_{2, i_1}+\ldots + e_{i_s}B_{2,i_s}, G)$ 
}
\end{prop}

\subsubsection{$D_{2n}$-covers}

We here introduce notation for dihedral covers.
Let $D_{2n}$ be the dihedral group of order $2n$. 
 In order to present $D_{2n}$, we use
 the notation
 \[
 D_{2n} = \langle \sigma, \tau \mid \sigma^2 = \tau^n = (\sigma\tau)^2 = 1\rangle.
 \]
 Given a $D_{2n}$-cover, we obtain a double cover, $D(X/Y)$, of $Y$ canonically by considering the
 $\CC(X)^{\tau}$-normalization of $Y$, where $\CC(X)^{\tau}$ denotes the fixed field 
 of the subgroup generated by $\tau$.  $X$ is an $n$-fold cyclic cover of $D(X/Y)$ and
 we denote these covering morphisms by
 $\beta_1(\pi) : D(X/Y) \to Y$ and $\beta_2(\pi) : X \to D(X/Y)$, respectively.

 \subsubsection{Elliptic surfaces and $D_{2p}$-covers of $\Sigma_d$}
 
 Let $\varphi : S \to \PP^1$ be an elliptic surface. As for the double cover diagram for $S$, we have 
 the following:
 
 \begin{itemize}
 
  \item $W$ is the Hirzebruch surface, 
 $\Sigma_d$, of degree 
 $d = 2\chi({\mathcal O}_S)$. Hence $d$ is always even. We denote $\widehat {W}$ by $\widehat{\Sigma}_d$.
 
 \item  $\Delta_{f'}$ is of the form $\Delta_0 + T$, where
 $\Delta_0$ is a section with $\Delta_0^2 = -d$ and $T \sim 3(\Delta_0 + d{\mathfrak f})$, 
 ${\mathfrak f}$ being a fiber of the ruling $\Sigma_d \to {\mathbb P}^1$.
 
 \end{itemize}

 In previous articles, \cite{tokunaga04, tokunaga12}, we studied $p$-cyclic covers ($p$: odd prime) 
 $g : X \to S$ such that $f\circ g : X \to \widehat{\Sigma}_d$ gives rise to a $D_{2p}$-cover of
  $\widehat{\Sigma}_{d}$. One of  the main results is as follows:
  
  \begin{thm}\label{thm:elliptic-dihed}{\cite[Theorem 3.1]{tokunaga12} Let $C$ be a reduced divisor on
  $S$ such that
  \begin{itemize}
   \item all irreducible components $C_i$ of $C = \sum_i C_i$ are horizontal
   with respect to the elliptic fibration (i.e., any irreducible component is not contained in any fiber), and 
%  \item If we denote $D = \sum_iD_i$ be irreducible decomposition of $D$.  
  \item $C$ and $\sigma_f^*C$ have no common component.
  \end{itemize}
 
  Then there exists a $p$-cyclic cover $g : X \to S$ such that
  \begin{enumerate}
  \item[(i)] $\Delta_g = C + \sigma_f^*C + \Xi + \sigma_f^*\Xi$, where $\Xi$ is effective
  divisor on $S$ such that  irreducible components of $\Xi$
 are all vertical and there is no common component between $\Xi$ and
    $\sigma_f^*\Xi$.
  \item[(ii)] $f\circ g : X \to \widehat{\Sigma}_d$ is a $D_{2p}$-cover of $\widehat{\Sigma}_d$
  such that $D(X/\widehat{\Sigma}_d) = S$,  $\beta_1(\pi)= f$ and $\beta_2(\pi) = g$,
  \end{enumerate}
  
  if and only if the following condition holds:
  
  Let $s(C_i) = \tilde{\psi}(C_i)$ $(i = 1,\ldots, r)$. There exist integers $a_i$ $(i = 1,\ldots, r)$ such that
\begin{itemize}
  \item $1 \le a_i < p$ $(i = 1, \ldots, r)$ and
  \item $\sum_{i=1}^r [a_i]s(C_i)$ is $p$-divisible in $\MW(S)$, i.e., 
  \[
  \sum_{i=1}^r [a_i]s(C_i) \in [p]\MW(S) := \{[p]s \mid s \in \MW(S)\}.
  \]
\end{itemize}  
}
\end{thm}

\begin{cor}\label{cor:main}{Suppose that $p\not| \MW(S)_{\tor}$ and $C$ consists of two components
$C_1$ and $C_2$. Then there exists a $p$-cyclic cover as in Theorem~\ref{thm:elliptic-dihed}
if and only the images of $s(C_1)$ and $s(C_2)$ are linearly dependent in $\MW(S)\otimes\ZZ/p\ZZ$.
}
\end{cor}

\begin{cor}\label{cor:main2}{Suppose that $p\not| \MW(S)_{\tor}$ and $C$ consists of two components
$C_1$ and $C_2$. Then there exists a $D_{2p}$-cover of $\widehat{\Sigma}_d$ branched at
$2\Delta_{f} + f(C)$
if and only the images of $s(C_1)$ and $s(C_2)$ are linearly dependent in $\MW(S)\otimes\ZZ/p\ZZ$.
}
\end{cor}

 %%%%%%%%%%%%%%%%%%%%%%%%%%%%%%%%%%
 
 %%%%%%%%%%%%%%%%%%%%%%%%%%%%%%%%%%

%Proof of Theorem~\ref{prescription}

\section{Proof of Theorem~\ref{prescription}}

We keep our notation from the Introduction.
Note that each surface $S_{p_i}$ has a singular fiber of type  $\I_2$ whose components arise from the exceptional divisor of the first blow up in $\mu_i$ which meets $O_{i}$ and the strict transform of the tangent line $l_{p_i}$ of $\mcQ$ at $p_i$. We will denote these components by $\Theta_{p_i,0}$ and $\Theta_{p_i,1}$.  All the other reducible singular fibers arise from the exceptional sets of the resolution $\bar{S}\rightarrow S''$, hence they are in 1 to 1 correspondence with the singularities of $\mcQ$. We will denote their components by $\Theta_{v,i}$ where $v\in \Sing(\mcQ)$. We will use the same symbol $\Theta_{v,i}$ for these components on both $S_{p_i}$.

 Let $D^\prime=s_1(D_1)+\cdots+s_2(D_n)$. Note that  the sum taken here is regarded as a sum of divisors on $S_{p_1}$. Then since the Abel-Jacobi map $\tilde{\psi}$ is a homomorphism, $\tilde{\psi}_1(D^\prime)=s_1(D_1)\dot+\cdots\dot+s_1(D_n)=C_1$. Hence by Lemma~\ref{lem:fund-relation} we have the equivalence
\[
D^\prime\underset{S_{p_1}}\sim (C_1)+d(O_1)+nF_1+a\Theta_{p_1,1}+\sum_{v\in \Sing(\mcQ)}\sum_{i=1}^{m_v-1} b_{v,i}\Theta_{v,i}
\]
Where $\underset{S_{p_1}}\sim$ denotes linear equivalence of divisors on $S_{p_1}$.

Similarly for $D^{\prime\prime}=s_2(D_1)+\cdots+s_2(D_n)$, we have the equivalence
\[
D^{\prime\prime}\underset{S_{p_2}}\sim (C_2)+d(O_2)+nF_2+a\Theta_{p_2,1}+\sum_{v\in \Sing(\mcQ)}\sum_{i=1}^{m_v-1} b_{v,i}\Theta_{v,i}.
\]
Note that by construction of $D^\prime, D^{\prime\prime}$, the coefficients $d, n, a, b_{v,i}$ are the same in both cases. Since ${\mu_{1\ast}}(D^\prime)=D_1+\cdots+D_m={\mu_{2\ast}}(D^{\prime\prime})$ and since linear equivalence behaves well under push-forward, we have 
\begin{align*}
\overline{C}_1+&nF_1+a\Theta_{p_1,1}+\sum_{v\in \Sing(\mcQ)}\sum_{i=1}^{m_v-1} b_{v,i}\Theta_{v,i}\\
&\underset{\bar{S}}{\sim} \overline{C}_2+nF_2+a\Theta_{p_2,1}+\sum_{v\in \Sing(\mcQ)}\sum_{i=1}^{m_v-1} b_{v,i}\Theta_{v,i}
\end{align*}
where we denote the images of $F_i$ and $\Theta_{p_i,1}$ by the same symbols.
Then since $F_1\underset{\bar{S}}{\sim}F_2$ and $\Theta_{p_1,1}\underset{\bar{S}}{\sim}\Theta_{p_2,1}$, because they are inverse images of lines of $\PP^2$, we obtain the equivalence
\[
\overline{C}_1\underset{\bar{S}}{\sim}\overline{C}_2.
\]
By pulling this equivalence back by ${\mu_1}$, we obtain
\begin{align*}
\widehat{C}_2&\underset{S_{p_1}}\sim C_1+\alpha O_1+\beta\Theta_{p_1,0}\\
&\underset{S_{p_1}}\sim C_1+\alpha O_1+\beta F-\beta \Theta_{p_1,1},
\end{align*}
for some integers $\alpha$ and $\beta$. Hence by Theorem~\ref{thm:shioda-basic} we have $\tilde{\psi}(\widehat{C}_2)=C_1$.
\qed

 %%%%%%%%%%%%%%%%%%%%%%%%%%%%%%%%%%
 
   %%%%%%%%%%%%%%%%%%%%%%%%%%%%%%%%%%%

%Construction of weak $k$-NT arrangement

%Construction of weak $k$-NT arrangement

\section{Construction for weak Namba-Tsuchihashi arrangements of type $k$ and $D_{2p}$-covers}

Let $\mcQ$ be a quartic consisting of $2$ irreducible conics intersecting at $4$ distinct points.
Let $f'' : S''\to \PP^2$ be a double cover of $\PP^2$ with $\Delta_{f''} = \mcQ$ and let 
\[
\begin{CD}
S'' @<{\mu_o}<< \bar{S} \\
@V{f''}VV                 @VV{f_o}V \\
\PP^2@<<{q_o}< \widehat{\PP}^2.
\end{CD}
\]
be the diagram for the canonical resolution of $S''$. Note that $q_o : \widehat{\PP}^2 \to \PP^2$ is 
a composition of blowing ups at the $4$ nodes $\Sing(\mcQ) = \{P_0, P_1, P_2, P_3\}$. 
We put $\tilde{f}_o = \mu_o\circ f'' = q_o\circ f_o$.
We choose a
point $x \in \mcQ$ such that 

{\bf Assumption.}
\begin{enumerate}

 \item[(i)] the tangent line $l_x$ at $x$ meets $\mcQ$ at two other distinct points and
 
 \item[(ii)] no tangent line at $\Sing(\mcQ)$ passes through $x$.
 
 \end{enumerate}
 
 Note that these conditions imply that $x$ is on a component of $\mathcal{Q}$ having degree greater than or equal to 2.
 
 The pencil of lines through $x$ gives rise to a pencil of curves of genus $1$, $\Lambda_x$, on $\bar{S}$.
 By resolving the base points,  $\mu_x : S_x \to \bar{S}$,  of $\Lambda_x$, we have a rational elliptic surface $\varphi_x : S_x \to \PP^1$.
 
 By our choice of $x$, $\varphi_x$ has $5$ singular fibers of type $\I_2$. Let $L_1, L_2$ and $L_3$
 be lines through $\{P_0, P_1\}, \{P_0, P_2\}$ and $\{P_0, P_3\}$. For each $L_i$,  $\mu_o^*L_i$ is of 
 the form $L_i^+ + L_i^-$. By our construction, $s_{L_i} := \mu_x^*L_i^+$ is  a section with
 $\langle s_{L_i}, s_{L_i} \rangle = 1/2$ for each $i$. Thus, by labeling the singular fibers suitably, we
 may assume that $s_{L_i}$ ($i = 1, 2, 3$) meet the singular fibers as follows:
 
 \bigskip
 
 \begin{center}
 %%%%%%%%%%%%%%%%%%%%%%%%%%%%

%WinTpicVersion4.23
\unitlength 0.1in
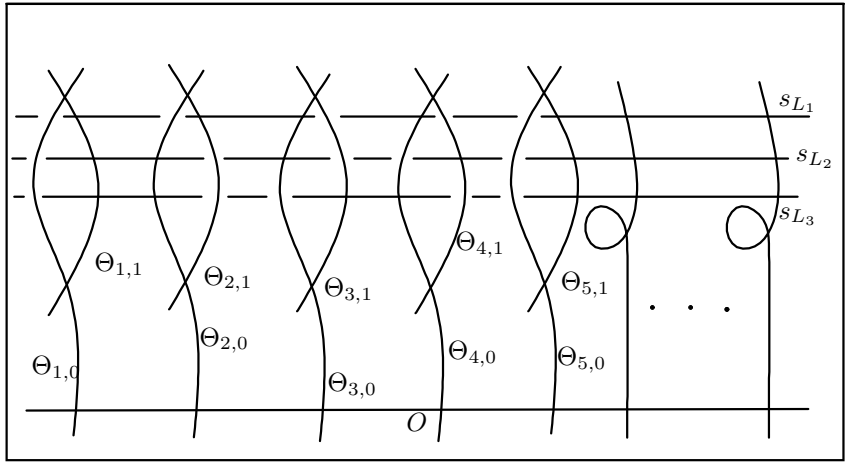
\begin{picture}( 43.8000, 23.9000)(  9.8000,-28.1000)
% BOX 1 0 3 0 Black White
% 2 980 420 5360 2810
% 
{{%
\special{pn 13}%
\special{pa 980 420}%
\special{pa 5360 420}%
\special{pa 5360 2810}%
\special{pa 980 2810}%
\special{pa 980 420}%
\special{pa 5360 420}%
\special{fp}%
}}%
% LINE 1 0 3 0 Black White
% 2 1084 2550 5175 2542
% 
{{%
\special{pn 13}%
\special{pa 1084 2550}%
\special{pa 5176 2542}%
\special{fp}%
}}%
% SPLINE 1 0 3 0 Black White
% 11 4180 830 4268 1244 4268 1497 4140 1696 4068 1696 4012 1558 4076 1481 4196 1543 4228 1688 4228 2661 4228 2692
% 
{{%
\special{pn 13}%
\special{pa 4180 830}%
\special{pa 4196 894}%
\special{pa 4204 924}%
\special{pa 4212 956}%
\special{pa 4218 988}%
\special{pa 4226 1018}%
\special{pa 4234 1050}%
\special{pa 4240 1082}%
\special{pa 4246 1112}%
\special{pa 4258 1176}%
\special{pa 4262 1206}%
\special{pa 4268 1238}%
\special{pa 4272 1268}%
\special{pa 4276 1300}%
\special{pa 4278 1332}%
\special{pa 4280 1362}%
\special{pa 4280 1394}%
\special{pa 4278 1426}%
\special{pa 4270 1494}%
\special{pa 4262 1528}%
\special{pa 4252 1560}%
\special{pa 4240 1594}%
\special{pa 4224 1624}%
\special{pa 4206 1650}%
\special{pa 4182 1674}%
\special{pa 4156 1690}%
\special{pa 4126 1700}%
\special{pa 4092 1702}%
\special{pa 4062 1694}%
\special{pa 4038 1674}%
\special{pa 4020 1644}%
\special{pa 4010 1610}%
\special{pa 4010 1574}%
\special{pa 4018 1540}%
\special{pa 4036 1510}%
\special{pa 4060 1488}%
\special{pa 4090 1480}%
\special{pa 4120 1484}%
\special{pa 4152 1500}%
\special{pa 4180 1522}%
\special{pa 4200 1550}%
\special{pa 4214 1580}%
\special{pa 4222 1610}%
\special{pa 4226 1642}%
\special{pa 4228 1674}%
\special{pa 4228 1708}%
\special{pa 4230 1740}%
\special{pa 4230 2476}%
\special{pa 4228 2508}%
\special{pa 4228 2692}%
\special{fp}%
}}%
% SPLINE 1 0 3 0 Black White
% 11 4920 830 5009 1244 5009 1497 4880 1696 4808 1696 4752 1558 4816 1481 4936 1543 4969 1688 4969 2661 4969 2692
% 
{{%
\special{pn 13}%
\special{pa 4920 830}%
\special{pa 4936 894}%
\special{pa 4944 924}%
\special{pa 4960 988}%
\special{pa 4966 1018}%
\special{pa 4974 1050}%
\special{pa 4980 1082}%
\special{pa 4986 1112}%
\special{pa 4992 1144}%
\special{pa 4998 1174}%
\special{pa 5004 1206}%
\special{pa 5008 1238}%
\special{pa 5012 1268}%
\special{pa 5016 1300}%
\special{pa 5018 1330}%
\special{pa 5020 1362}%
\special{pa 5020 1394}%
\special{pa 5018 1426}%
\special{pa 5016 1460}%
\special{pa 5010 1492}%
\special{pa 5004 1526}%
\special{pa 4994 1560}%
\special{pa 4980 1594}%
\special{pa 4964 1624}%
\special{pa 4946 1650}%
\special{pa 4924 1672}%
\special{pa 4896 1690}%
\special{pa 4866 1700}%
\special{pa 4834 1702}%
\special{pa 4802 1694}%
\special{pa 4778 1674}%
\special{pa 4760 1646}%
\special{pa 4750 1610}%
\special{pa 4750 1574}%
\special{pa 4758 1540}%
\special{pa 4776 1510}%
\special{pa 4800 1488}%
\special{pa 4828 1480}%
\special{pa 4860 1484}%
\special{pa 4890 1500}%
\special{pa 4918 1522}%
\special{pa 4940 1550}%
\special{pa 4954 1578}%
\special{pa 4964 1610}%
\special{pa 4968 1642}%
\special{pa 4970 1674}%
\special{pa 4970 1738}%
\special{pa 4972 1772}%
\special{pa 4972 2348}%
\special{pa 4970 2380}%
\special{pa 4970 2692}%
\special{fp}%
}}%
% STR 2 0 3 0 Black White
% 4 5020 914 5020 990 2 0 0 0
% $s_{L_1}$
\put(50.2000,-9.9000){\makebox(0,0)[lb]{$s_{L_1}$}}%
% STR 2 0 3 0 Black White
% 4 5110 1204 5110 1280 2 0 0 0
% $s_{L_2}$
\put(51.1000,-12.8000){\makebox(0,0)[lb]{$s_{L_2}$}}%
% STR 2 0 3 0 Black White
% 4 5020 1503 5020 1580 2 0 0 0
% $s_{L_3}$
\put(50.2000,-15.8000){\makebox(0,0)[lb]{$s_{L_3}$}}%
% STR 2 0 3 0 Black White
% 4 3074 2580 3074 2657 2 0 0 0
% $O$
\put(30.7400,-26.5700){\makebox(0,0)[lb]{$O$}}%
% STR 2 0 3 0 Black White
% 4 1110 2324 1110 2400 2 0 0 0
% $\Theta_{1,0}$
\put(11.1000,-24.0000){\makebox(0,0)[lb]{$\Theta_{1,0}$}}%
% STR 2 0 3 0 Black White
% 4 3330 1674 3330 1750 2 0 0 0
% $\Theta_{4,1}$
\put(33.3000,-17.5000){\makebox(0,0)[lb]{$\Theta_{4,1}$}}%
% STR 2 0 3 0 Black White
% 4 1440 1784 1440 1860 2 0 0 0
% $\Theta_{1,1}$
\put(14.4000,-18.6000){\makebox(0,0)[lb]{$\Theta_{1,1}$}}%
% STR 2 0 3 0 Black White
% 4 3290 2244 3290 2320 2 0 0 0
% $\Theta_{4,0}$
\put(32.9000,-23.2000){\makebox(0,0)[lb]{$\Theta_{4,0}$}}%
% SPLINE 1 0 3 0 Black White
% 5 2010 770 1750 1380 1930 1920 1960 2690 1960 2690
% 
{{%
\special{pn 13}%
\special{pa 2010 770}%
\special{pa 1992 800}%
\special{pa 1972 830}%
\special{pa 1952 858}%
\special{pa 1934 888}%
\special{pa 1916 916}%
\special{pa 1880 976}%
\special{pa 1862 1004}%
\special{pa 1846 1034}%
\special{pa 1818 1094}%
\special{pa 1804 1122}%
\special{pa 1792 1152}%
\special{pa 1782 1182}%
\special{pa 1772 1210}%
\special{pa 1764 1240}%
\special{pa 1758 1270}%
\special{pa 1754 1300}%
\special{pa 1750 1360}%
\special{pa 1750 1388}%
\special{pa 1758 1448}%
\special{pa 1764 1478}%
\special{pa 1772 1508}%
\special{pa 1802 1598}%
\special{pa 1814 1628}%
\special{pa 1826 1660}%
\special{pa 1840 1690}%
\special{pa 1852 1720}%
\special{pa 1866 1750}%
\special{pa 1878 1782}%
\special{pa 1890 1812}%
\special{pa 1904 1842}%
\special{pa 1914 1874}%
\special{pa 1926 1904}%
\special{pa 1936 1936}%
\special{pa 1944 1968}%
\special{pa 1952 1998}%
\special{pa 1964 2062}%
\special{pa 1970 2092}%
\special{pa 1978 2156}%
\special{pa 1984 2252}%
\special{pa 1984 2348}%
\special{pa 1982 2380}%
\special{pa 1982 2412}%
\special{pa 1980 2444}%
\special{pa 1978 2478}%
\special{pa 1972 2574}%
\special{pa 1968 2606}%
\special{pa 1966 2638}%
\special{pa 1962 2672}%
\special{pa 1960 2690}%
\special{fp}%
}}%
% SPLINE 1 0 3 0 Black White
% 3 1200 770 1460 1380 1200 2050
% 
{{%
\special{pn 13}%
\special{pa 1200 770}%
\special{pa 1218 800}%
\special{pa 1238 830}%
\special{pa 1256 858}%
\special{pa 1272 888}%
\special{pa 1290 918}%
\special{pa 1308 946}%
\special{pa 1340 1006}%
\special{pa 1356 1034}%
\special{pa 1384 1094}%
\special{pa 1398 1122}%
\special{pa 1410 1152}%
\special{pa 1430 1212}%
\special{pa 1438 1240}%
\special{pa 1446 1270}%
\special{pa 1452 1300}%
\special{pa 1460 1360}%
\special{pa 1460 1418}%
\special{pa 1456 1478}%
\special{pa 1444 1538}%
\special{pa 1436 1566}%
\special{pa 1428 1596}%
\special{pa 1408 1656}%
\special{pa 1396 1686}%
\special{pa 1382 1716}%
\special{pa 1370 1746}%
\special{pa 1354 1776}%
\special{pa 1340 1806}%
\special{pa 1292 1896}%
\special{pa 1256 1956}%
\special{pa 1240 1986}%
\special{pa 1204 2046}%
\special{pa 1200 2050}%
\special{fp}%
}}%
% SPLINE 1 0 3 0 Black White
% 5 1380 760 1120 1370 1300 1910 1330 2680 1330 2680
% 
{{%
\special{pn 13}%
\special{pa 1380 760}%
\special{pa 1362 790}%
\special{pa 1342 820}%
\special{pa 1322 848}%
\special{pa 1304 878}%
\special{pa 1286 906}%
\special{pa 1250 966}%
\special{pa 1232 994}%
\special{pa 1216 1024}%
\special{pa 1188 1084}%
\special{pa 1174 1112}%
\special{pa 1162 1142}%
\special{pa 1152 1172}%
\special{pa 1142 1200}%
\special{pa 1134 1230}%
\special{pa 1128 1260}%
\special{pa 1124 1290}%
\special{pa 1120 1350}%
\special{pa 1120 1378}%
\special{pa 1128 1438}%
\special{pa 1134 1468}%
\special{pa 1142 1498}%
\special{pa 1172 1588}%
\special{pa 1184 1618}%
\special{pa 1196 1650}%
\special{pa 1210 1680}%
\special{pa 1222 1710}%
\special{pa 1236 1740}%
\special{pa 1248 1772}%
\special{pa 1260 1802}%
\special{pa 1274 1832}%
\special{pa 1284 1864}%
\special{pa 1296 1894}%
\special{pa 1306 1926}%
\special{pa 1314 1958}%
\special{pa 1322 1988}%
\special{pa 1334 2052}%
\special{pa 1340 2082}%
\special{pa 1348 2146}%
\special{pa 1354 2242}%
\special{pa 1354 2338}%
\special{pa 1352 2370}%
\special{pa 1352 2402}%
\special{pa 1350 2434}%
\special{pa 1348 2468}%
\special{pa 1342 2564}%
\special{pa 1338 2596}%
\special{pa 1336 2628}%
\special{pa 1332 2662}%
\special{pa 1330 2680}%
\special{fp}%
}}%
% SPLINE 1 0 3 0 Black White
% 3 1830 740 2090 1350 1830 2020
% 
{{%
\special{pn 13}%
\special{pa 1830 740}%
\special{pa 1848 770}%
\special{pa 1868 800}%
\special{pa 1886 828}%
\special{pa 1902 858}%
\special{pa 1920 888}%
\special{pa 1938 916}%
\special{pa 1970 976}%
\special{pa 1986 1004}%
\special{pa 2014 1064}%
\special{pa 2028 1092}%
\special{pa 2040 1122}%
\special{pa 2060 1182}%
\special{pa 2068 1210}%
\special{pa 2076 1240}%
\special{pa 2082 1270}%
\special{pa 2090 1330}%
\special{pa 2090 1388}%
\special{pa 2086 1448}%
\special{pa 2074 1508}%
\special{pa 2066 1536}%
\special{pa 2058 1566}%
\special{pa 2038 1626}%
\special{pa 2026 1656}%
\special{pa 2012 1686}%
\special{pa 2000 1716}%
\special{pa 1984 1746}%
\special{pa 1970 1776}%
\special{pa 1922 1866}%
\special{pa 1886 1926}%
\special{pa 1870 1956}%
\special{pa 1834 2016}%
\special{pa 1830 2020}%
\special{fp}%
}}%
% SPLINE 1 0 3 0 Black White
% 5 2670 790 2410 1400 2590 1940 2620 2710 2620 2710
% 
{{%
\special{pn 13}%
\special{pa 2670 790}%
\special{pa 2652 820}%
\special{pa 2632 850}%
\special{pa 2612 878}%
\special{pa 2594 908}%
\special{pa 2576 936}%
\special{pa 2540 996}%
\special{pa 2522 1024}%
\special{pa 2506 1054}%
\special{pa 2478 1114}%
\special{pa 2464 1142}%
\special{pa 2452 1172}%
\special{pa 2442 1202}%
\special{pa 2432 1230}%
\special{pa 2424 1260}%
\special{pa 2418 1290}%
\special{pa 2414 1320}%
\special{pa 2410 1380}%
\special{pa 2410 1408}%
\special{pa 2418 1468}%
\special{pa 2424 1498}%
\special{pa 2432 1528}%
\special{pa 2462 1618}%
\special{pa 2474 1648}%
\special{pa 2486 1680}%
\special{pa 2500 1710}%
\special{pa 2512 1740}%
\special{pa 2526 1770}%
\special{pa 2538 1802}%
\special{pa 2550 1832}%
\special{pa 2564 1862}%
\special{pa 2574 1894}%
\special{pa 2586 1924}%
\special{pa 2596 1956}%
\special{pa 2604 1988}%
\special{pa 2612 2018}%
\special{pa 2624 2082}%
\special{pa 2630 2112}%
\special{pa 2638 2176}%
\special{pa 2644 2272}%
\special{pa 2644 2368}%
\special{pa 2642 2400}%
\special{pa 2642 2432}%
\special{pa 2640 2464}%
\special{pa 2638 2498}%
\special{pa 2632 2594}%
\special{pa 2628 2626}%
\special{pa 2626 2658}%
\special{pa 2622 2692}%
\special{pa 2620 2710}%
\special{fp}%
}}%
% SPLINE 1 0 3 0 Black White
% 3 2500 750 2760 1360 2500 2030
% 
{{%
\special{pn 13}%
\special{pa 2500 750}%
\special{pa 2518 780}%
\special{pa 2538 810}%
\special{pa 2556 838}%
\special{pa 2572 868}%
\special{pa 2590 898}%
\special{pa 2608 926}%
\special{pa 2640 986}%
\special{pa 2656 1014}%
\special{pa 2684 1074}%
\special{pa 2698 1102}%
\special{pa 2710 1132}%
\special{pa 2730 1192}%
\special{pa 2738 1220}%
\special{pa 2746 1250}%
\special{pa 2752 1280}%
\special{pa 2760 1340}%
\special{pa 2760 1398}%
\special{pa 2756 1458}%
\special{pa 2744 1518}%
\special{pa 2736 1546}%
\special{pa 2728 1576}%
\special{pa 2708 1636}%
\special{pa 2696 1666}%
\special{pa 2682 1696}%
\special{pa 2670 1726}%
\special{pa 2654 1756}%
\special{pa 2640 1786}%
\special{pa 2592 1876}%
\special{pa 2556 1936}%
\special{pa 2540 1966}%
\special{pa 2504 2026}%
\special{pa 2500 2030}%
\special{fp}%
}}%
% SPLINE 1 0 3 0 Black White
% 5 3290 790 3030 1400 3210 1940 3240 2710 3240 2710
% 
{{%
\special{pn 13}%
\special{pa 3290 790}%
\special{pa 3272 820}%
\special{pa 3252 850}%
\special{pa 3232 878}%
\special{pa 3214 908}%
\special{pa 3196 936}%
\special{pa 3160 996}%
\special{pa 3142 1024}%
\special{pa 3126 1054}%
\special{pa 3098 1114}%
\special{pa 3084 1142}%
\special{pa 3072 1172}%
\special{pa 3062 1202}%
\special{pa 3052 1230}%
\special{pa 3044 1260}%
\special{pa 3038 1290}%
\special{pa 3034 1320}%
\special{pa 3030 1380}%
\special{pa 3030 1408}%
\special{pa 3038 1468}%
\special{pa 3044 1498}%
\special{pa 3052 1528}%
\special{pa 3082 1618}%
\special{pa 3094 1648}%
\special{pa 3106 1680}%
\special{pa 3120 1710}%
\special{pa 3132 1740}%
\special{pa 3146 1770}%
\special{pa 3158 1802}%
\special{pa 3170 1832}%
\special{pa 3184 1862}%
\special{pa 3194 1894}%
\special{pa 3206 1924}%
\special{pa 3216 1956}%
\special{pa 3224 1988}%
\special{pa 3232 2018}%
\special{pa 3244 2082}%
\special{pa 3250 2112}%
\special{pa 3258 2176}%
\special{pa 3264 2272}%
\special{pa 3264 2368}%
\special{pa 3262 2400}%
\special{pa 3262 2432}%
\special{pa 3260 2464}%
\special{pa 3258 2498}%
\special{pa 3252 2594}%
\special{pa 3248 2626}%
\special{pa 3246 2658}%
\special{pa 3242 2692}%
\special{pa 3240 2710}%
\special{fp}%
}}%
% SPLINE 1 0 3 0 Black White
% 3 3120 760 3380 1370 3120 2040
% 
{{%
\special{pn 13}%
\special{pa 3120 760}%
\special{pa 3138 790}%
\special{pa 3158 820}%
\special{pa 3176 848}%
\special{pa 3192 878}%
\special{pa 3210 908}%
\special{pa 3228 936}%
\special{pa 3260 996}%
\special{pa 3276 1024}%
\special{pa 3304 1084}%
\special{pa 3318 1112}%
\special{pa 3330 1142}%
\special{pa 3350 1202}%
\special{pa 3358 1230}%
\special{pa 3366 1260}%
\special{pa 3372 1290}%
\special{pa 3380 1350}%
\special{pa 3380 1408}%
\special{pa 3376 1468}%
\special{pa 3364 1528}%
\special{pa 3356 1556}%
\special{pa 3348 1586}%
\special{pa 3328 1646}%
\special{pa 3316 1676}%
\special{pa 3302 1706}%
\special{pa 3290 1736}%
\special{pa 3274 1766}%
\special{pa 3260 1796}%
\special{pa 3212 1886}%
\special{pa 3176 1946}%
\special{pa 3160 1976}%
\special{pa 3124 2036}%
\special{pa 3120 2040}%
\special{fp}%
}}%
% STR 2 0 3 0 Black White
% 4 2010 1854 2010 1930 2 0 0 0
% $\Theta_{2,1}$
\put(20.1000,-19.3000){\makebox(0,0)[lb]{$\Theta_{2,1}$}}%
% STR 2 0 3 0 Black White
% 4 1990 2174 1990 2250 2 0 0 0
% $\Theta_{2,0}$
\put(19.9000,-22.5000){\makebox(0,0)[lb]{$\Theta_{2,0}$}}%
% STR 2 0 3 0 Black White
% 4 2660 2414 2660 2490 2 0 0 0
% $\Theta_{3,0}$
\put(26.6000,-24.9000){\makebox(0,0)[lb]{$\Theta_{3,0}$}}%
% STR 2 0 3 0 Black White
% 4 2640 1914 2640 1990 2 0 0 0
% $\Theta_{3,1}$
\put(26.4000,-19.9000){\makebox(0,0)[lb]{$\Theta_{3,1}$}}%
% DOT 1 0 3 0 Black White
% 4 4360 2010 4560 2010 4750 2020 4750 2020
% 
{{%
\special{pn 4}%
\special{sh 1}%
\special{ar 4360 2010 10 10 0  6.28318530717959E+0000}%
\special{sh 1}%
\special{ar 4560 2010 10 10 0  6.28318530717959E+0000}%
\special{sh 1}%
\special{ar 4750 2020 10 10 0  6.28318530717959E+0000}%
\special{sh 1}%
\special{ar 4750 2020 10 10 0  6.28318530717959E+0000}%
}}%
% SPLINE 1 0 3 0 Black White
% 3 3710 750 3970 1360 3710 2030
% 
{{%
\special{pn 13}%
\special{pa 3710 750}%
\special{pa 3728 780}%
\special{pa 3748 810}%
\special{pa 3766 838}%
\special{pa 3782 868}%
\special{pa 3800 898}%
\special{pa 3818 926}%
\special{pa 3850 986}%
\special{pa 3866 1014}%
\special{pa 3894 1074}%
\special{pa 3908 1102}%
\special{pa 3920 1132}%
\special{pa 3940 1192}%
\special{pa 3948 1220}%
\special{pa 3956 1250}%
\special{pa 3962 1280}%
\special{pa 3970 1340}%
\special{pa 3970 1398}%
\special{pa 3966 1458}%
\special{pa 3954 1518}%
\special{pa 3946 1546}%
\special{pa 3938 1576}%
\special{pa 3918 1636}%
\special{pa 3906 1666}%
\special{pa 3892 1696}%
\special{pa 3880 1726}%
\special{pa 3864 1756}%
\special{pa 3850 1786}%
\special{pa 3802 1876}%
\special{pa 3766 1936}%
\special{pa 3750 1966}%
\special{pa 3714 2026}%
\special{pa 3710 2030}%
\special{fp}%
}}%
% SPLINE 1 0 3 0 Black White
% 5 3880 750 3620 1360 3800 1900 3830 2670 3830 2670
% 
{{%
\special{pn 13}%
\special{pa 3880 750}%
\special{pa 3862 780}%
\special{pa 3842 810}%
\special{pa 3822 838}%
\special{pa 3804 868}%
\special{pa 3786 896}%
\special{pa 3750 956}%
\special{pa 3732 984}%
\special{pa 3716 1014}%
\special{pa 3688 1074}%
\special{pa 3674 1102}%
\special{pa 3662 1132}%
\special{pa 3652 1162}%
\special{pa 3642 1190}%
\special{pa 3634 1220}%
\special{pa 3628 1250}%
\special{pa 3624 1280}%
\special{pa 3620 1340}%
\special{pa 3620 1368}%
\special{pa 3628 1428}%
\special{pa 3634 1458}%
\special{pa 3642 1488}%
\special{pa 3672 1578}%
\special{pa 3684 1608}%
\special{pa 3696 1640}%
\special{pa 3710 1670}%
\special{pa 3722 1700}%
\special{pa 3736 1730}%
\special{pa 3748 1762}%
\special{pa 3760 1792}%
\special{pa 3774 1822}%
\special{pa 3784 1854}%
\special{pa 3796 1884}%
\special{pa 3806 1916}%
\special{pa 3814 1948}%
\special{pa 3822 1978}%
\special{pa 3834 2042}%
\special{pa 3840 2072}%
\special{pa 3848 2136}%
\special{pa 3854 2232}%
\special{pa 3854 2328}%
\special{pa 3852 2360}%
\special{pa 3852 2392}%
\special{pa 3850 2424}%
\special{pa 3848 2458}%
\special{pa 3842 2554}%
\special{pa 3838 2586}%
\special{pa 3836 2618}%
\special{pa 3832 2652}%
\special{pa 3830 2670}%
\special{fp}%
}}%
% STR 2 0 3 0 Black White
% 4 3880 1884 3880 1960 2 0 0 0
% $\Theta_{5,1}$
\put(38.8000,-19.6000){\makebox(0,0)[lb]{$\Theta_{5,1}$}}%
% STR 2 0 3 0 Black White
% 4 3860 2274 3860 2350 2 0 0 0
% $\Theta_{5,0}$
\put(38.6000,-23.5000){\makebox(0,0)[lb]{$\Theta_{5,0}$}}%
% LINE 1 0 3 0 Black White
% 12 1030 1010 1140 1010 1280 1010 1780 1010 1920 1010 2590 1010 2730 1010 3220 1010 3320 1010 3650 1010 3780 1010 5180 1010
% 
{{%
\special{pn 13}%
\special{pa 1030 1010}%
\special{pa 1140 1010}%
\special{fp}%
\special{pa 1280 1010}%
\special{pa 1780 1010}%
\special{fp}%
\special{pa 1920 1010}%
\special{pa 2590 1010}%
\special{fp}%
\special{pa 2730 1010}%
\special{pa 3220 1010}%
\special{fp}%
\special{pa 3320 1010}%
\special{pa 3650 1010}%
\special{fp}%
\special{pa 3780 1010}%
\special{pa 5180 1010}%
\special{fp}%
}}%
% LINE 1 0 3 0 Black White
% 12 1010 1230 1080 1230 1190 1230 2010 1230 2120 1230 2670 1230 2800 1230 2990 1230 3110 1230 3570 1230 3700 1230 5070 1230
% 
{{%
\special{pn 13}%
\special{pa 1010 1230}%
\special{pa 1080 1230}%
\special{fp}%
\special{pa 1190 1230}%
\special{pa 2010 1230}%
\special{fp}%
\special{pa 2120 1230}%
\special{pa 2670 1230}%
\special{fp}%
\special{pa 2800 1230}%
\special{pa 2990 1230}%
\special{fp}%
\special{pa 3110 1230}%
\special{pa 3570 1230}%
\special{fp}%
\special{pa 3700 1230}%
\special{pa 5070 1230}%
\special{fp}%
}}%
% LINE 1 0 3 0 Black White
% 12 1020 1430 1070 1430 1170 1430 2000 1430 2140 1430 2340 1430 2480 1430 3300 1430 3420 1430 3560 1430 3690 1430 5120 1430
% 
{{%
\special{pn 13}%
\special{pa 1020 1430}%
\special{pa 1070 1430}%
\special{fp}%
\special{pa 1170 1430}%
\special{pa 2000 1430}%
\special{fp}%
\special{pa 2140 1430}%
\special{pa 2340 1430}%
\special{fp}%
\special{pa 2480 1430}%
\special{pa 3300 1430}%
\special{fp}%
\special{pa 3420 1430}%
\special{pa 3560 1430}%
\special{fp}%
\special{pa 3690 1430}%
\special{pa 5120 1430}%
\special{fp}%
}}%
\end{picture}%

 %%%%%%%%%%%%%%%%%%%%%%%%%%%%
 \end{center}
 
 \bigskip
 
Here $\Theta_{1,1}$ is the component arising from $l_x$. 
 Put $s_i := [2]s_{L_i}$. By \cite[Theorem 9.1]{kodaira}, we infer that $s_iO = 0$ ($i = 1, 2, 3$). Hence
 we have
 \[
 \begin{array}{ccc}
 \langle s_i, s_i \rangle = 2 \quad (i = 1, 2, 3) & \quad\quad & \langle s_i, s_j \rangle = 0 \quad (1 \le i < j \le 3).
 \end{array}
 \]
 Hence $\MW(S)^0 ( \cong A_1^{\oplus 3}) = \ZZ s_1\oplus \ZZ s_2 \oplus \ZZ s_3$, where
 $\MW(S)^0$ denote the narrow part of $\MW(S)$.

 \begin{lem}\label{lem:3-conics}{ Let $\mcQ$ and $x \in \mcQ$ be as above.  There exists
 exactly $3$ conics $\mcC_{1, x}, \mcC_{2, x}$ and $\mcC_{3, x}$ through $x$ such 
 that (i) $\mcC_{i, x}$ does 
 not pass through $P_0, P_1, P_2$ and $P_3$ for each $i$ and (ii) for each 
 $z_o \in \mcC_{i, x}\cap \mcQ$, the intersection multiplicity at $z_o$, $I_{z_o}(C_{i, x}, \mcQ)$ is
 even.
 }
 \end{lem}
 
 \proof Put $\mcC_{i, x} := \tilde {f}_o\circ\mu_x(s_i) =  \tilde {f}_o\circ\mu_x([-1]s_i)$ 
 ($i = 1, 2, 3$). Since $s_i \in \MW(S)^0$, 
 $s_i$ always meets the identity component of a fiber for each $i$.  For each $i$,
  $\mcC_{i, x}$ passes $x$ and meets a general line through $x$ at another point distinct from $x$. Hence
  $\mcC_{i, x}$ is a smooth conic. Since $s_i \neq [-1]s_i$, for each 
 $z_o \in \mcC_{i, x}\cap \mcQ$, the intersection multiplicity at $z_o$, 
 $I_{z_o}(\mcC_{i, x}, \mcQ)$ is even.

Conversely, any conic $\mcC$ satisfying the conditions gives rise to sections $s_{\mcC}, 
\, [-1]s_{\mcC} \in \MW(S)^0$ with $\langle s_{\mcC}, s_{\mcC} \rangle = \langle [-1]s_{\mcC}, [-1]s_{\mcC} \rangle = 2$. Thus $\mcC$ is one of $\mcC_{i,x}$ ($i = 1, 2, 3$). \qed

\begin{lem}\label{lem:key-1}{Choose $x, y \in \mcQ$ satisfying the Assumption. 
%Put
%$\MW(S_x)^0 = \ZZ s_{1,x}\oplus \ZZ s_{2,x} \oplus \ZZ s_{3, x}$ and
%$\MW(S_y)^0 = \ZZ s_{1, y}\oplus \ZZ s_{2, y} \oplus \ZZ s_{3, y}$, where
Put $s_{i, x} = [2]s_{L_i, x}$ and $s_{i, y} = [2]s_{L_i, y}$ $(i = 1, 2, 3)$, where
$s_{L_i, x}$ and $s_{L_i, y}$ are the sections arising from $(\tilde {f}_o\circ \mu_x)^*(L_i)$ and
$(\tilde {f}_o\circ \mu_y)^*(L_i)$, respectively. 
Put $\mcC_{i, x} = \tilde{f}\circ\mu_x(s_{i,x})$ and $\mcC_{i, y} = \tilde{f}\circ\mu_x(s_{i, y})$
Let $\mcC_{i,x}^{\pm}$ (resp. $\mcC_{i,y}^{\pm}$) be bisections on $S_y$ 
(resp. $S_x$) arising $\mcC_{i, x}$ (resp. $\mcC_{i,y}$).
Then $s(C_{i, y}^+) = s_{i, x}$ or $s(C_{i, y}^+) = [-1]s_{i, x}$ on $\MW(S_x)$ and
$s(C_{i, x}^+) = s_{i, y}$ or $s(C_{i, x}^+) = [-1]s_{i, y}$ on $\MW(S_y)$.
}
\end{lem}
 
 \proof Our statement is immediate by Theorem~\ref{prescription}. \qed

\bigskip

Fix a positive integer $k \ge 2$ and take
 $(k_1, k_2, k_3) \in \ZZ_{\ge 0}^{\oplus 3}$ with $k_1 \ge k_2 \ge k_3$ and 
$k = k_1 + k_2 + k_3$. Choose $x_i$ ($i = 1, \ldots, k$) $\in \mcQ$ satisfying
the Assumption. We consider $k$ conics $\mcC_1, \ldots, \mcC_k$ as follows:
\[
\mcC_i := \left \{ \begin{array}{cc}
                  \tilde{f}_o\circ\mu_{x_i}(s_{1, x_i}) & i = 1, \ldots, k_1 \\
                  \tilde{f}_o\circ\mu_{x_i}(s_{2, x_i}) & i = k_1 + 1, \ldots, k_1 + k_2 \\
                  \tilde{f}_o\circ\mu_{x_i}(s_{3, x_i}) & i = k_1 + k_2 + 1, \ldots, k.
                  \end{array}
                  \right .
 \]
 Here we choose $s_{j, x_i}$ $(j = 1, 2, 3)$ as in Lemma~\ref{lem:key-1}. Then
 $\mcQ + \sum_{i=1}^k\mcC_i$ is a weak $k$-Namba-Tsuchihashi arrangement.   
  \bigskip
  
  \begin{lem}\label{lem:key-2}{Let $p$ be an odd prime. For $\mcC_i$ and $\mcC_j$ as
  above, there exists a $D_{2p}$-cover of $\PP^2$ branched at $2\mcQ + p(\mcC_i + 
  \mcC_j)$ if and only if $i, j \in \{1, \ldots, k_1\}$, $\{k_1 + 1, \ldots, k_1 + k_2\}$ or
  $\{k_1 + k_2 + 1, \ldots, k\}$.
  }
  \end{lem}
  
  \proof Choose a general point $x \in \mcQ\setminus\{x_1, \ldots, x_k\}$ satisfying the
  Assumption. Let $\varphi_x : S_x \to \PP^1$ be a rational elliptic surface as before and 
  let
 \[
\begin{CD}
S_{x}' @<{\mu^\prime_x}<< S_x \\
@V{f'_x}VV                 @VV{f_x}V \\
\Sigma_2@<<{q^\prime_x}< \widehat{\Sigma}_2.
\end{CD}
\]
be its double cover diagram. By our construction of $S_x$, we infer that
 there exists a composition of blowing ups 
 $q_x : \widehat{\Sigma}_2 \to \widehat{\PP}^2$ such that
 \[
\begin{CD}
S'' @<{\mu_o}<< \bar{S} @<{\mu_x}<< S_x \\
@V{f'_o}VV                 @VV{f_o}V @VV{f_x}V \\
\PP^2@<<{q_o}< \widehat{\PP}^2 @<<{q_x}< \widehat{\Sigma}_2.
\end{CD}
\]

Let $\overline{\mcC}_i$ be the strict transform of $\mcC_i$ with respect to $q_o\circ q_x$.
By considering the Stein factorization, there exists a $D_{2p}$-cover of $\PP^2$
branched at $2\mcQ + p(\mcC_i + \mcC_j)$ if and only if  
there exists a $D_{2p}$-cover of $\widehat{\Sigma}_2$
branched at $2\Delta_{f_x} + p(\overline{\mcC}_i + \overline{\mcC}_j)$.
Put $f_x^*\overline{\mcC}_i = \overline{\mcC}_i^+ + \overline{\mcC}_i^-$. By our choice of
$\mcC_i$,  we have          
\[
s(\mcC_j^+) = \left \{ \begin{array}{cc}
                                \mbox{$s_{1, x}$ or $[-1]s_{1, x}$} & \mbox{if $1\le j \le k_1$} \\
                                \mbox{$s_{2, x}$ or $[-1]s_{2, x}$} & \mbox{if $k_1 + 1\le j \le k_1 + k_2$} \\
                                \mbox{$s_{3, x}$ or $[-1]s_{3, x}$} & \mbox{if $k_1 + k_2 + 1\le j \le k$} 
                                \end{array} \right .
\]
Hence by Corollary~\ref{cor:main2}, our statement follows.
\qed

Let $>_{\lex}$ be the {\it lexicographic order} on $\ZZ_{\ge 0}^{\oplus 3}$ (see \cite{CLO}). Take
$(k_1, k_2, k_3), (k'_1, k'_2, k'_3) \in \ZZ_{\ge 0}^{\oplus 3}$ as in Lemma~\ref{lem:key-2} such that
$(k_1, k_2, k_3) >_{\lex} (k'_1, k'_2, k'_3)$. Choose 

\begin{itemize}

\item $\mcC_1^{(1)}, \ldots, \mcC_k^{(1)}$ for $(k_1, k_2, k_3)$, and

\item  $\mcC_1^{(2)}, \ldots, \mcC_k^{(2)}$ for $(k'_1, k'_2, k'_3)$

\end{itemize}
as above. Put $\mcB_1 = \mcQ + \sum_{i=1}^k\mcC_i^{(1)}$ and $\mcB_2 = \mcQ + \sum_{i=1}^k
\mcC_i^{(2)}$.  Note that both $\mcB_1$ and $\mcB_2$ are weak $k$-Namba-Tsuchihashi arrangements. 

\begin{prop}\label{prop:main-wNT}{There exists no homeomorphism $h: (\PP^2, \mcB_1) \to (\PP^2, 
\mcB_2)$ such that $h(\mcQ) = \mcQ$.
}
\end{prop}

\proof Suppose that a homeomorphism $h : (\PP^2, \mcB_1) \to (\PP^2, \mcB_2)$ exists. 
Put
\[
\begin{array}{ccc}
h(\mcC_1^{(1)}) = \mcC_{i_1}^{(2)}, &  & h(\mcC_{k_1+1}^{(1)}) = \mcC_{i_2}^{(2)}.
\end{array}
\]
Since $\Cov_b(\PP^2, 2\mcQ + p(\mcC_1^{(1)} +\mcC_j^{(1)}), D_{2p}) \neq \emptyset$  $(j = 2, \ldots, k_1)$, 
by
Proposition~\ref{prop:1-1}, there exists at least $(k_1 -1)$ $C_j^{(2)}$ different form $C_{i_1}^{(2)}$ such that
$\Cov_b(\PP^2, 2\mcQ + p(\mcC_{i_1}^{(2)} + \mcC_j^{(2)}), D_{2p}) \neq \emptyset$. As
$k_1 \ge k'_1$, we infer that $k_1 =k'_1$ by Lemma~\ref{lem:key-2}. 
Similarly, since $\Cov_b(\PP^2, 2\mcQ + p(\mcC_{k_1+1}^{(1)} +\mcC_j^{(1)}), D_{2p}) \neq \emptyset$ 
 $(j = k_1 + 1, \ldots, k_1 + k_2)$, there exists at least $(k_2 -1)$ $C_j^{(2)}$ different form 
 $C_{i_2}^{(2)}$ such that
$\Cov_b(\PP^2, 2\mcQ + p(\mcC_{i_2}^{(2)} + \mcC_j^{(2)}), D_{2p}) \neq \emptyset$ and we infer
that $k_2 = k'_2$. This contradicts  $(k_1, k_2, k_3) >_{\lex} (k'_1, k'_2, k'_3)$. \qed

\begin{cor}\label{cor:main-wNT}{If both $\mcB_1$ and $\mcB_2$ are $k$-NT arrangements for $k \ge 3$,
then $(\mcB_1, \mcB_2)$ is a Zariski pair.}
\end{cor}
\proof If both $\mcB_1$ and $\mcB_2$ are $k$-NT arrangements and $k \ge 3$, $h(\mcQ) =\mcQ$ holds 
for any homeomorphism $(\PP^2, \mcB_1) \to (\PP^2, \mcB_2)$. \qed.

%%%%%%%%%%%%%%%%%%%%%%%%%%%%%%%%%%

   %%%%%%%%%%%%%%%%%%%%%%%%%%%%%%%%

 \section{Proof of Theorem~\ref{thm:zariski-kplet}}

 %'±'̈ꕶ'ð‰Á'¦'é21 may 
In this section we will construct $y(k, 3)$ Namba-Tsuchihashi arrangements of type $k$ which form
 Zariski $y(k, 3)$-plets for conic arrangements.  The main ingredient of the proof is an explicit method to calculate the equations of the conics that appeared in Proposition~\ref{prop:main-wNT}.   
 We keep the notation in \S 3.%families of

\subsection{Explicit method in finding equations of bisections}\label{explicit}
Let $S$ be an elliptic surface over $\PP^1$, whose generic fiber is given by a
Weirstrass form  $y^2=f(t,x)$.  Note that the degree of $f$ with respect to $x$, $\deg_xf(t,x)$, is
equal to $3$.
Consider a rational point $P$ of the generic fiber of $S$ with coordinates $(x(t), y(t))$. Let $L$ be the line passing through $P$ in $\mathbb{A}_{\CC(t)}^2$ defined by
\[
L: y=r(t)(x-x(t))+y(t)
\]
for some $r(t)\in \CC(t)$. Then, since $(x(t), y(t))$ is a rational point of $S$ and since $L$ passes through $(x(t), y(t))$, the equation
\[
\left\{r(t)(x- x(t))+ y(t))\right\}^2-f(t,x)
\]
factors into the form
\[
\left\{r(t)(x-x(t))+y(t))\right\}^2-f(x,t)=(x-x(t))g_{}(t,x).
\]
Since $\deg_xg_{}(t,x)=2$, the support of the intersection of $L$ and $g_{}(t,x)=0$, viewed as rational functions of $S$ over the field $\CC$, defines  a bisection $D$ on $S$, and by the definition of the Abel-Jacobi map, $s(D)=-P$. By varying $r(t)$, we get a family of bisections $\mathcal{D}$ such that any $D\in\mathcal{D}$ has the same image $s(D)=-P$ under the Abel-Jacobi map.

In general, $r(t)$ can be any rational function in $\CC(t)$ and hence $\mathcal{D}$ becomes an enormously large family. We will find various useful strata of $\mathcal{D}$ by restricting $r(t)$.

Note that the images  of the bisections, obtained by the method  above, in $\PP^2$ is defined by $g(t,x)=0$.

\subsection{Construction of $y(k,3)$ Namba-Tsuchihashi arrangements of type $k$}
Let $[T:X:Z]$ be homogeneous coordinates of $\PP^2$ and let $t=T/Z$ and $x=X/Z$. Let $\mcQ: (XZ-T^2+2Z^2) (X^2-2XZ+T^2-4Z^2)=0$, and  $z_o=[0:1:0]$. We will consider the  elliptic surface $S$, as in \S 3, branched along $\mcQ$ with blow-up center $z_o=[0:1:0]$. The Weierstrass equation of $S$ is
\[
y^2=(x-t^2+2)(x^2-2x+t^2-4).
\]
%$S$ is a rational elliptic surface whose configuration of singular fibers is $5I_2$, $2I_1$. 
%By \cite{oguiso-shioda} the Mordell-Weil lattice of $S$ is isomorphic to
% $(A_1^\ast)^{\oplus3}\oplus\ZZ/2\ZZ$. 
The singularities of $\mcQ$ consists of four nodes at $\{[\pm 1, -1, 1], [\pm 2, 2, 1]\}$ and we put $P_0 = [2, 2, 1], P_1=[-1, -1, 1], 
P_2=[1, -1, 1], P_3=[-2, 2, 1]$. Then the lines $L_1, L_2$ and $L_3$ 
in \S 3 are given by
  $L_1: X-T=0$, $L_2:X-3T+4Z=0$ and $L_3: X-2Z=0$, respectively.
    These lines gives rise to sections $s_{L_1}$, $s_{L_2}$ and
  $s_{L_3}$ as in \S 3. 
%   Each of  $l_1, l_2$ and $l_3$ passes through $[2,2,1]$ and one other intersection point of $C_1$
%  and $C_2$, hence give rise to sections $ s_{l_1}^{\pm}, s_{l_2}^\pm, s_{l_3}^\pm$ in $S$. Let 
%  $s_1=s_{l_1}^+$, $s_2=s_{l_2}^+$ and $s_3=s_{l_3}^+$.  We
%may assume that the intersection of  $s_1, s_2$ and $s_3$ with the singular fibers of $S$ is as 
%in the following figure:
%
%\begin{center}
%\input No24-2
%\end{center}
%
%By the explicit formula for the height pairing, we find that 
%$\langle s_i,s_i\rangle=\frac{1}{2}$ $(i=1,2,3)$, $\langle s_i,s_j\rangle=0$ $(i\not=j)$. Hence we can
% assume that $s_1, s_2, s_3$ is a generator of the free part of $\MW(S)$. 
Also the conic $XZ-T^2+2Z^2=0$ gives rise to the two-torsion section, which we will denote by $s_t$. These
sections are given explicitly as follows:
\begin{align*}
s_{L_1}&=\left(t,\sqrt{-2}(t+1)(t-2)\right)\\
s_{L_2}&=\left(3t-4, \sqrt{-10}(t-1)(t-2)\right)\\
s_{L_3} &=\left(2,\sqrt{-1}(t+2)(t-2)\right)\\
s_t&=\left(t^2+2,0\right)
\end{align*}   

Put $s_i = [2]s_{L_i}$.  
Explicitly, $s_1$, $s_2$ and $s_3$ are given as follows:
\begin{align*}
s_1&=\left(\frac{1}{2}\,{t}^{2}-2,-\frac{\sqrt {-2}}{4} t \left( {t}^{2}-4 \right) \right)\\
s_2&=\left(\frac{1}{10}\,{t}^{2}-2,-{\frac {3\sqrt{-10}}{100}} t \left( {t}^{2}+20\right)\right)\\
s_3&=\left({t}^{2}-{\frac {17}{4}},-\frac{3\sqrt{-1}}{8} \left( 4\,{t}^{2}-19 \right) \right)
\end{align*}

As we see in \S 3,  each $s_i$ gives rise to a conic in $\PP^2$ that is tangent to 
$\mcQ$ at $z_o$ and three other points.
The defining equation of the conic $\mcC_1$ (resp. $\mcC_2$, $\mcC_3$) corresponding to   $s_1$ (resp. $s_2$, $s_3$) is    $\mcC_1: XZ-\frac{1}{2}T^2+2Z^2=0$ (resp.  $\mcC_2: XZ-\frac{1}{10}T^2+2Z^2=0$, $\mcC_3:XZ-T^2+2Z^2=0$).

By applying the method of finding the explicit equations of bisections 
whose image of the Abel-Jacobi map is  $[-1]s_1$ in
\ref{explicit},   we obtain a family of bisections 
$\bar{\mcD}_1$ whose images in $\PP^2$ under $\tilde{f}_o\circ \mu_{z_o}$ are curves given by  $\bar{g}_{1,r}(t,x)=0$, $ r\in\CC(t)$
where \begin{align*}
\bar{g}_{1,r}(t,x)=&\frac{\,{t}^{4}}{4}+\frac{\,r\sqrt {-2}}{2}\,{t}^{3}+\frac{1}{2}\,{t}^{2}x- \left( 3+\frac{\,{r}}{2}
^{2} \right) {t}^{2}\\&\quad-{x}^{2}-2\,r\sqrt {-2}\,t+ \left( {r}^{2}+2
 \right) x+2\,{r}^{2}+4
\end{align*}
By specializing $r$ to $r=\frac{1}{2}\sqrt{-2} t-a_1$  $(a_1\in \CC)$, we obtain a sub-family $\mathcal{D}_1\subset\mathcal{\bar{D}}_1$ of bisections. The defining equations of the images under $\bar{q}\circ f$ specialize to  $g_{1,a_1}(t,x)$, $a_1\in \CC$ where $g_{1,a_1}(t,x)$ is given by
\[
g_{1,a_1}(t,x)=\left( -2-\frac{a_1^2}{2} \right) {t}^{2}-xt\sqrt {-2}a_1-{x}^{2}+ \left( 
{a_1}^{2}+2 \right) x+2\,{a_1}^{2}+4.
\]
By construction, any $D_{1,a}\in \mathcal{D}_1$ satisfies $s(D_{a_1})=s_1$. 

Similarly, by applying the same method to $[-1]s_2$ and $[-1]s_3$, we obtain families of bisections $\mathcal{D}_{2}$ (resp. $\mathcal{D}_{3}$) parametrized by $a_2\in\CC$ (resp. $a_3\in\CC$) such that any $D_{i,a_i}\in\mathcal{D}_{i}$ satisfies $s(D_{i,a_i})=s_i$ $(i = 2, 3)$ and the defining equations of their images in $\PP^2$ are given by
\begin{align*}
g_{2,a_2}(t,x)&=\frac 1{10}\,{a_2}^{2}{t}^{2}-10\,{t}^{2}-{\frac {12}{5}}\,\sqrt{-10}ta_2+2\,{a_2}
^{2}- \frac 35\,\sqrt{-10}xta_2+{a_2}^{2}x+4+2\,x-{x}^{2},
\\
g_{3,a_3}(t,x)&=-{a_3}^{2}{t}^{2}+{\frac {17}{4}}\,{a_3}^{2}+{a_3}^{2}x-3\,\sqrt{-1}a_3{t}^{2}+{\frac 
{57}{4}}\,\sqrt{-1} a_3+ \frac 54\,{t}^{2}-{\frac {161}{16}}+{\frac {17}{4}}\,x-{x}^{2}
\end{align*}
We will denote the image of $D_{i,a_i}$  in $\PP^2$ by  $\mcC_{i,a_i}$, which is defined by $g_{i,a_i}(t,x)=0$. 
\begin{lem}
Under the notation above,  the following statements hold:
\begin{enumerate}
\item For $i=1, 2, 3$ and general $a_i$, $\mcC_{i,a_i}$ is a smooth conic. In particular there are only a finite number of non-reduced members of $\mathcal{D}_i$.
 \item For $i=1, 2, 3$, there exist only a finite number of  $\mcC_{i,a_i}$ passing through any given point $p\in\PP^2$.
 \item For $i=1, 2, 3$, there exist only a finite number of $a_i$ such that $\mcC_{i,a_i}$ and $\mcQ$ have a intersection point with multiplicity $\ge 4$. 
  \item For all $\{i,j\}\subset\{1,2,3\}$ and any given $a_i\in\CC$ such that $\mcC_{i,a_i}$ is reduced, there  exist only a finite number of $a_j$  such that $\mcC_{j,a_j}$ and  $\mcC_{i,a_i}$ do not intersect transversally.
\end{enumerate}
\end{lem}
\proof Each statement can be proved by direct calculations. The computer algebra software MAPLE was used in the calculations.

\begin{lem}\label{strong-kNT}
Let $k$ be an integer $\geq 3$ and $k_1, k_2, k_3$ be non-negative integers such that $k=k_1+k_2+k_3$. Then there exists $\alpha_1,\cdots,\alpha_k\in\CC$ such that the configuration 
\[\mcQ+\sum_{i=1}^{k_1}\mcC_{1,\alpha_i}+\sum_{i=k_1+1}^{k_1+k_2}\mcC_{2,\alpha_i}+\sum_{i=k_1+k_2+1}^{k}\mcC_{3,\alpha_i}\]
 becomes a Namba-Tsuchihashi configuration of type $k$ satisfying the properties described before Proposition~\ref{prop:main-wNT}.
\end{lem}

\proof This follows directly from the previous lemma. 

By combining Theorem \ref{prescription}, Lemma \ref{strong-kNT} and Proposition~\ref{prop:main-wNT}, we obtain 
Theorem \ref{thm:zariski-kplet}.

\noindent Shinzo BANNAI and  Hiro-o TOKUNAGA\\
Department of Mathematics and Information Sciences\\
Graduate School of Science and Engineering,\\
Tokyo Metropolitan University\\
1-1 Minami-Ohsawa, Hachiohji 192-0397 JAPAN \\
{\tt shinzo.bannai@gmail.com},
{\tt tokunaga@tmu.ac.jp}
      
  \end{document}